\documentclass[proceedings]{stacs}
\stacsheading{2009}{337--348}{Freiburg}
\firstpageno{337}

\usepackage{amssymb}
\usepackage{latexsym}
\usepackage{amsmath}
\usepackage{amsfonts}
\usepackage{mathrsfs}

\theoremstyle{plain}
\theoremstyle{definition}
\theoremstyle{question}\newtheorem{question}[thm]{Question}

\newcommand{\rst}[1]{\ensuremath{{\mathbin\mid}%
\raise-.5ex\hbox{$#1$}}}  

\newcommand{\n}{\mathbb{N}}
\newcommand{\z}{\mathbb{Z}}
\newcommand{\az}{A^{\mathbb{Z}}}
\newcommand{\bz}{B^{\mathbb{Z}}}
\newcommand{\cz}{C^{\mathbb{Z}}}
\newcommand{\AZ}{(\az,F)}
\newcommand{\BZ}{(\bz,G)}
\newcommand{\CZ}{(\cz,H)}

\begin{document}

\title[Undecidable properties of limit set dynamics of Cellular Automata]{Undecidable properties of limit set dynamics of Cellular Automata}

\author[lab1]{P. Di Lena}{Pietro Di Lena}
\address[lab1]{Department of Computer Science,
University of Bologna,
Mura Anteo Zamboni 7,
Bologna, Italy
  \newline }  
\email[P. Di Lena]{dilena@cs.unibo.it}  
\email[L. Margara]{margara@cs.unibo.it} 

\author[lab1]{L. Margara}{Luciano Margara}


\subjclass{Theory of Computation, Computation by Abstract Devices}


\begin{abstract}
  \noindent Cellular Automata (CA) are discrete dynamical systems and an abstract model of parallel computation. The limit set of a cellular automaton is its maximal topological attractor. A well know result, due to Kari, says that all nontrivial properties of limit sets are undecidable. In this paper we consider properties of limit set dynamics, i.e. properties of the dynamics of Cellular Automata restricted to their limit sets. There can be no equivalent of Kari's Theorem for limit set dynamics. Anyway we show that there is a large class of undecidable properties of limit set dynamics, namely all properties of limit set dynamics which imply stability or the existence of a unique subshift attractor. As a consequence we have that it is undecidable whether the cellular automaton map restricted to the limit set is the identity, closing, injective, expansive, positively expansive, transitive.\end{abstract}

\maketitle

\vspace*{-5mm}
\section*{Introduction}\label{S:one}

\emph{Cellular Automata} (CA) are discrete dynamical systems and, at the same time, an abstract model of parallel computation. Every cellular automaton has a finite description in terms of a finite block mapping called \emph{local rule}. A general problem for CA is to determine what are the properties which are algorithmically decidable/undecidable given the local rule.

The limit set $\Omega_F$ of a cellular automaton $\AZ$ is the set of all configurations which occur after arbitrarily long iterates of the CA map, i.e. $x\in\Omega_F$ if and only if $\forall n\in \n, F^{-n}(x)\neq \emptyset$. The limit set is the maximal topological attractor of a cellular automaton (then it is always nonempty and closed) and it is fundamental to understand the long-term behavior of such systems.
Kari's Theorem \cite{Kari94} says that all nontrivial properties of limit sets are undecidable. This implies, for example, that we cannot decide algorithmically if some given configuration is in the limit set or not and we cannot even decide if some given word is contained in some configuration of the limit set.
Kari's undecidability theorem  uniquely regards properties of the configurations contained in the limit set, but it does not include properties of the dynamics of Cellular Automata restricted to their limit set. The motivation of this work is to try to understand what are the undecidable properties of the limit set dynamics, i.e. properties of the dynamical systems $(\Omega_F,F)$.
It is easy to find simple examples of nontrivial decidable properties of $F:\Omega_F\rightarrow \Omega_F$ which imply that Kari's Theorem cannot be extended to whole limit set dynamics.
Anyway, we can show that there is a large and interesting class of properties of $F:\Omega_F\rightarrow \Omega_F$ which are undecidable. 
For instance, we show that any property of limit set dynamics which implies stability or the existence of a unique subshift attractor is undecidable. Stated in another way, we obtain that any decidable property of limit set dynamics must be a property of some unstable cellular automaton with at least two subshift attractors. As a consequence we show that it is not possible to decide algorithmically whether the 
cellular automaton map restricted to the limit set is the identity, closing, injective, expansive, positively expansive and transitive.

The paper is organized as follows. In Section \ref{preliminaries} we provide the basic background in Symbolic Dynamics and Cellular Automata needed to understand the rest of the paper.
In Section \ref{properties} we formally define what properties of limit sets are and we show some preliminary results. In Section \ref{results} we discuss our main results. Section \ref{conclusions} is devoted to concluding remarks.

\section{Preliminaries}\label{preliminaries}

\subsection{Symbolic Dynamics}
In this section we review only those notions which are strictly necessary to understand our proofs.
See  \cite{LM95} for a complete introduction to Symbolic Dynamics.

Let $A$ be a finite alphabet with at least two elements. We denote by $A^n$ the set of words 
of length $n$ over $A$,  by $A^*=\cup_{n\in\n} A^n$ the set of words over $A$ and by $\az$ the set of doubly infinite sequences $(x_{i})_{i\in \z}$ of symbols $x_{i} \in A$. We denote by $x_{[i,j]}\in A^{j-i+1}$ the subword $x_{i}x_{i+1}...x_{j}$. We use the shortcut $w \sqsubset x$ to say that $w\in A^{+}$ is a subword of $x \in \az$. 

Define a \emph{metric} $d$ on $\az$ by $d(x,y) = 2^{-n}$ where $n = min\{|i| \mid x_{i}\neq y_{i}\}$. The set $\az$ endowed with metric $d$ is a compact metric space. For $u \in A^{*}$ and $i \in \z$, denote by $[u]_{i} = \{x\in \az \mid x_{[i,i+|u|-1]} = u\}$
a \emph{cylinder set}. For a lighter notation, we will refer to the cylinder set $[u]_{i}$ simply by $[u]$. A cylinder set is a clopen (closed and open) set in $\az$. Every clopen set in $\az$ is a finite union of cylinder sets. 

The \emph{shift map} $\sigma : \az \rightarrow \az$ is defined by $\sigma(x)_{i} = x_{i+1}$. 
The shift map is continuous and biiective on $\az$. The dynamical system $(\az, \sigma)$ is
called \emph{full shift}. A \emph{shift space} or \emph{subshift} is a non-empty closed subset
$\Sigma \subseteq \az$ which is strongly shift invariant, i.e. $\sigma(\Sigma)=\Sigma$.
We will usually denote the \emph{shift dynamical system} $(\Sigma, \sigma)$
simply with $\Sigma$. A subshift $\Sigma$ is a \emph{zero-dimensional space}, i.e. for every two
different points $x,y\in\Sigma$ there exists disjoint clopen sets $U,V \subset \Sigma$ such that $x\in U, y\in V$.

We denote by $\mathcal{L}_n(\Sigma)=\{w\in A^n \mid  \exists x \in \Sigma, w\sqsubset x\}$ the set of words of length $n$ of the subshift $\Sigma$. The \emph{language} of $\Sigma$ is defined by
$\mathcal{L}(\Sigma) = \cup_{n\in\n}\mathcal{L}_n(\Sigma)$. Any subshift $\Sigma$ is completely determined by the set of its \emph{forbidden words} $A^{*}\setminus\mathcal{L}(\Sigma)$.  A \emph{shift of finite type} (SFT) is a subshift which can be defined by a \emph{finite} set of forbidden words. Let $\Sigma$ be a subshift on alphabet $A$. We denote by $\Sigma_k=\{x\in\az \mid \forall i\in\z, x_{[i,i+k)}\in\mathcal{L}_k(\Sigma)\}$ the \emph{SFT approximation of order $k>0$} of $\Sigma$. Note that $\forall k>0, \Sigma\subseteq \Sigma_k$ and that $\Sigma_k$ is a SFT since it is defined by the
finite set of forbidden words $A^k\setminus \mathcal{L}_k(\Sigma)$. If $\Sigma$ is a SFT then there exists some $k>0$ such that $\forall k'\geq k,\Sigma=\Sigma_{k'}$. We say that the least such $k>0$ is the \emph{order} of $\Sigma$. A generalization of SFTs are \emph{sofic shifts}. A subshift $S$ is sofic if and only if its language $\mathcal{L}(S)$ is \emph{regular}.
A subshift $\Sigma$ is \emph{mixing} if there exists $n>0$ such that for all clopen sets $U,V\subseteq \Sigma$, $\sigma^{n}(U)\cap V\neq \emptyset$. 

Let $\Sigma_1, \Sigma_2$ be subshifts. A \emph{factor map}  $F:\Sigma_1 \rightarrow \Sigma_2$ is
a continuous, onto, $\sigma$-commuting mapping. A factor map is actually a \emph{block code},
i.e. $F$ is induced by some \emph{$k$-block mapping} $f : \mathcal{L}^k_1(\Sigma_1) \rightarrow \mathcal{L}_1(\Sigma_2)$ where $k>0$. The mixing and sofic properties are preserved under factor maps.

A factor map $F$ is \emph{right-closing} if $x,y \in \Sigma, x_{(-\infty,i]}=y_{(-\infty,i]}$ and $F(x)=F(y)$ imply $x=y$. The definition of left-closing is equivalent. By using a simple compactness argument it is possible to prove that closing is equivalent to the following condition:
$\exists n> 0$ such that $\forall x,y \in \Sigma, \forall i \in \z$ if $x_{[i,i+n)}=y_{[i,i+n)}$ and  $F(x)_{[i,i+2n]}=F(y)_{[i,i+2n]}$ then $x_{i+n}=y_{i+n}$. The closing property imposes strong constraint on the mapping. For example, it is possible to prove that if $\Sigma$ is a mixing SFT and $F: \Sigma \rightarrow \Sigma$ is continuous, $\sigma$-commuting and closing then $F$ is onto, i.e. $F(\Sigma)=\Sigma$. 

An endomorphism $F: \Sigma \rightarrow \Sigma$ is \emph{positively expansive} if there exists $\epsilon > 0$ such that for all distinct $x,y\in\Sigma$ there exists $n\in \n$ such that $d(F^n(x),F^n(y))>\epsilon$.
If $F$ is invertible then it is \emph{expansive} if there exists $\epsilon > 0$
such that for all distinct $x,y\in\Sigma$ there exists $n\in \z$ such that $d(F^n(x),F^n(y))>\epsilon$.
Both expansive and positively expansive endomorphisms of subshifts must be closing.
The map $F$ is \emph{transitive}, if for any nonempty open sets $U,V\subseteq \Sigma$ there exists $n\in\n$ such that $F^{-n}(U)\cap V\neq\emptyset$. Both expansive and positively expansive endomorphisms of mixing SFT are transitive.

\subsection{Cellular Automata}

One-dimensional \emph{Cellular Automata} (CA) are endomorphisms of full shifts. We denote CA by pairs $\AZ$ where $F : \az \rightarrow \az$ is some continuous and $\sigma$-commuting function. The \emph{global rule} $F$ is a $(2r+1)$-block map, i.e. there exists some
\emph{local rule} $f : A^{2r+1} \rightarrow A$ of \emph{radius} $r \geq 0$ such that
\begin{center}
$\forall x \in \az, F(x)_{i} = f(x_{i-r},...,x_{i+r})$.
\end{center}
It is sometimes useful to extend the local rule to the finite-block mapping
\begin{center}
$f^*: A^k \rightarrow A^{k-2r}$ for every $k\geq 2r+1$,
\end{center}
such that
\begin{center}
$f^*(x_1,...,x_k)=f(x_{1},...,x_{2r+1})f(x_{2},...,x_{2r+2})..f(x_{k-2r},...,x_{k})$.
\end{center}

Our investigation regards properties of the \emph{limit behavior} of Cellular Automata. To understand the limit behavior the concept of \emph{attractor} is fundamental. An attractor is a nonempty closed set which attracts the orbits of its neighboring points. 

\begin{definition} Let $\AZ$ be a cellular automaton. The \emph{$\omega$-limit} of a set $U \subseteq \az$ with respect to $F$ is defined by
$\omega_F(U) = \cap_{n>0}\overline{\cup_{m>n} F^{m}(U)}$. 
\end{definition}
When it is clear from the context, we will denote the $\omega$-limit simply with $\omega$.
In zero-dimensional spaces the following two definitions of attractors are equivalent.

\begin{definition} Let $\AZ$ be a cellular automaton.  A nonempty closed set $Y\subseteq \az$ such that $F(Y)=Y$ is an \emph{attractor} of $\AZ$
\begin{itemize}
\item [$1.$] if $\forall \epsilon>0, \exists \delta>0$ such that $\forall x\in \az$
\begin{displaymath}
d(x,Y)<\delta \implies \forall n\in\n, d(F^n(x),Y)<\epsilon\;\;and\;\;\lim_{n\to \infty} d(F^n(x),Y)=0.
\end{displaymath}
\item [$2.$] if and only if $Y=\omega(U)$ where $U$ is a clopen $F$-invariant set, i.e.
$F(U)\subseteq U$.
\end{itemize}
\end{definition}
A useful property of attractors is that every neighborhood of an attractor contains a clopen $F$-invariant set whose $\omega$-limit is the attractor itself. We show the proof for completeness. We first need a general lemma.

\begin{lemma} \label{limit} Let $\AZ$ be a cellular automaton and let $U,V\subseteq \az$ be clopen sets. Assume that $\forall x\in U, \exists n_x\in \n$ such that $F^{n_x}(x)\in V$. Then there exists $n\in \n$ such that $\forall x\in U, \exists n_x\leq n, F^{n_x}(x)\in V$.
\end{lemma}
\proof{For $i\in\n$ define $X_i=\{x\in U \mid \forall j\leq i, F^j(x)\notin V\}$. Since $U,V$ are clopen it follows that for every  $i\in \n, X_i$ is clopen and $X_i \supseteq X_{i+1}$.
Assume that for every $i\in \n, X_i\neq\emptyset$ then, by compactness, $X=\cap_{i\in\n} X_i$ is nonempty which implies that there exists $x \in U$ such that $\forall i\in \n, F^i(x)\notin V$ contradicting
the hypothesis.\qed}

\begin{proposition} \label{attr} Let $\AZ$ be a cellular automaton and let  $Y\subseteq\az$ be an attractor. Then for every $\epsilon>0$ there is an $F$-invariant clopen set $U\subseteq \mathcal{B}_\epsilon(Y)$ such that $\omega(U)=Y$.
\end{proposition}
\proof{For $\epsilon>0$, denote $Y_\epsilon=\mathcal{B}_\epsilon(Y)$. Note that for every $\epsilon>0$, $Y_\epsilon$ is a clopen set.
Choose some $\epsilon>0$. By definition, there is some $0< \delta<\epsilon$ such that
\begin{displaymath}
x\in Y_\delta \implies \forall n, F^n(x)\in Y_\epsilon \;\;and\;\;  \lim_{n \to \infty} d(F^n(x),Y)=0. 
\end{displaymath}
Choose some $0<\epsilon_0<\delta$ then there is some $0<\delta_0<\epsilon_0$ such that
\begin{displaymath}
x\in Y_{\delta_0} \implies \forall n, F^n(x)\in Y_{\epsilon_0} \;\;and\;\;  \lim_{n \to \infty} d(F^n(x),Y)=0. 
\end{displaymath}

If $x\in Y_\delta$ then $\displaystyle\lim_{n \to \infty} d(F^n(x),Y)=0$ so there is some $n_x\in \n$ such that $F^{n_x}(x)\in Y_{\delta_0}$. By Lemma \ref{limit}, there is some $n\in\n$
such that for every $x\in Y_\delta, \exists n_x\leq n, F^{n_x}(x)\in Y_{\delta_0}$ then $\forall x \in Y_\delta, F^{n}(x)\in Y_{\epsilon_0}\subseteq Y_\delta$. We obtained that there is some $n\in \n$ such that $F^n(Y_\delta)\subseteq Y_\delta$ then $Y_\delta$ is $F^n$-invariant.
We now define a clopen set $U\subseteq Y_\epsilon$ which is $F$-invariant.

Let $x\in Y_\delta$, since $F^n(Y_\delta)\subseteq Y_\delta$ and $Y_\delta$ is clopen, there is a word $w \sqsubset x$ such that $[w]\subseteq Y_\delta$ and
$[(f^*)^n(w)] \subseteq Y_\delta$ (where the length of  $(f^*)^n(w)$ is greater than $0$). In particular, since $Y_\delta$ is the union of a finite collection of cylinders, there is a finite set of words $w^0_1, ..., w^0_{k_0}$ such that $Y_\delta = [w^0_1] \cup ... \cup  [w^0_{k_0}]$ and $[(f^*)^n(w^0_i)]\subseteq Y_\delta$ for $1\leq i \leq k_0$.
By considering iterates of $f^*$ on such words we can obtain a sequence of clopen sets $U_0, U_1, .., U_n$ such that $F(U_j)\subseteq U_{j+1}$. Set $U_0 = Y_\delta=[w^0_1] \cup ... \cup  [w^0_{k_0}]$ and define the clopen set $U_1=\cup^{k_0}_{i=1}([f^*(w^0_{i})]\cap Y_\epsilon)=[w^1_1] \cup ... \cup  [w^1_{k_1}]$. Note that for every $i \in [0,k_0]$ we have $F([w^0_i])\subseteq [f^*(w^0_i)], F([w^0_i])\subseteq Y_\epsilon$ and $[f^*(w^0_{i})]\cap Y_\epsilon$ is clopen. Then $F(U_0) \subseteq U_1 \subseteq Y_\epsilon$. Iterating for $j \in [1,n]$ we obtain the sequence of clopen sets
$U_j = \cup^{k_{j-1}}_{i=1}([f^*(w^{j-1}_{i})]\cap Y_\epsilon)=[w^j_1] \cup ... \cup  [w^j_{k_j}]$ such that $F(U_{j-1}) \subseteq U_j \subseteq Y_\epsilon$.
Now define $U=\cup^n_{j=0}U_j$. We have that $U\subseteq Y_\epsilon$ is clopen, $F(U)\subseteq U$ and $\omega(U)=Y$.\qed}

In the context of Cellular Automata, a particular class of attractors are those attractors which are also subshifts.

\begin{definition} Let $\AZ$ be a cellular automaton. A nonempty closed set $Y\subseteq \az$ is a \emph{subshift attractor} if it is an attractor and if $\sigma(Y)=Y$.
\end{definition}
The following two propositions characterize subshift attractors of CA. 

\begin{definition} Let $\AZ$ be a cellular automaton. We say that a clopen and $F$-invariant set $U\subseteq\az$ is \emph{spreading} if there exists some
$k>0$ such that $F^k(U)\subseteq \sigma^{-1}(U)\cap U\cap \sigma(U)$.
 \end{definition}

\begin{proposition} \label{spreading} \cite{FK07} Let $\AZ$ be a cellular automaton and let $U\subseteq \az$ be a clopen and $F$-invariant set. Then $\omega(U)$ is a subshift attractor if and only if $U$ is spreading.
\end{proposition}

\begin{proposition}  \label{subshift} \cite{FK07} Let $\AZ$ be a cellular automaton and let $U\subseteq \az$ be a clopen $F$-invariant spreading set. Then there exists a mixing SFT $\Sigma \subseteq U$ with the following properties:
\begin{itemize}
\item $F(\Sigma)\subseteq\Sigma$
\item $W=\{[w] \mid w \in \mathcal{L}_{k}(\Sigma), k$ is the order of $\Sigma\}$ is clopen and $F$-invariant
\item $\omega(\Sigma)=\omega(W)=\omega(U)$.
\end{itemize}
\end{proposition}
\noindent
A cellular automaton has at least one attractor which is called \emph{limit set}.

\begin{definition} Let $\AZ$ be a cellular automaton. The \emph{limit set} of $\AZ$ is defined by
$\Omega_F=\cap^{\infty}_{i=0}F^{i}(\az)=\omega_F(\az)$.
\end{definition}
Note that a configuration $x\in\az$ is in the limit set if and only if for every $n~\in~\n$, $F^{-n}(x)\neq\emptyset$. The limit set is a subshift attractor and it is also the \emph{maximal attractor}, i.e. every other attractor is contained in the limit set. The limit set can be the \emph{unique attractor}. In particular, if the map is transitive on the limit set then it is the unique attractor (the converse is not true). An attractor is  a \emph{minimal attractor} if it does not contain any proper subset which is an attractor. 
A unique attractor is both maximal and minimal. There is a very simple class of minimal attractors. We say that a state $s\in A$ is \emph{spreading} if the local rule has the property $f(x_1, ..., x_{2r+1})=s$ if $\exists x_i=s$. If a cellular automaton has a \emph{spreading state} $s$ then the clopen set $[s]$ is
$F$-invariant and spreading and $\omega([s])=\{...sss....\}$ is a minimal subshift attractor.

Cellular Automata limit sets received great attention. Here we review just some basic facts. One question which is still not well understood concerns the class of subshifts which can be limit sets of CA (see, for example, \cite{Hurd90, Maass95}). The most immediate distiction is between limit sets of stable and
unstable CA. A cellular automaton $\AZ$ is called \emph{stable} if there exists some $n\in\n$ such that
$F^n(\az)=\Omega_F$. It is called \emph{unstable} otherwise. The limit sets of stable CA are mixing sofic shifts since they are factors of full shifts. There are sofic subshifts which are limit sets of unstable CA but no limit set of unstable CA can be a SFT \cite{Hurd90}. It is actually unknown whether a subshift can be the limit set of both a stable and of an unstable CA. The simplest example of limit set subshift is the subshift consisting of just one configuration. A cellular automaton whose limit set is a single configuration is called
\emph{nilpotent}. 
If a cellular automaton is nilpotent then the unique configuration in the limit set must be fixed by both
$\sigma$ and $F$ then the automaton must be stable. 

\section{Properties of limit sets}\label{properties}

An important aspect of CA is that they can be enumerated. Every cellular automaton is described by its local rule. Local rules are defined by a finite amount of information and, in particular, for any fixed radius and cardinality of the alphabet there are only finitely many possible CA local rules. 

Choose some enumeration function for CA local rules. We denote by $\#\AZ\in\n$ the \emph{rule number} associated to $\AZ$. A \emph{property} $\mathcal{P}$ of CA is a collection of CA rule numbers. A property is called \emph{trivial} if either all CA have such property or none has.  A property $\mathcal{P}$ is \emph{decidable} whether there exists some algorithm such that, for any given $\AZ$, it
always computes if either $\#\AZ\in\mathcal{P}$ or $\#\AZ\notin\mathcal{P}$. A subclass of CA properties are, in particular, properties of the limit sets. 

\begin{definition} A property $\mathcal{P}$ is a \emph{property of limit sets} if and only if  the following condition holds: if
$\#\AZ\in\mathcal{P}$ and $\BZ$ is a cellular automaton such that $\Omega_F=\Omega_G$ then $\#\BZ\in\mathcal{P}$.
\end{definition}

Nilpotency is a property of limit sets. While a property of limit sets is always a property of CA, the converse is not always true. For example, surjectivity is not a property of the limit sets since it is easy to construct a surjective CA and a not surjective CA which have the same limit set. For example, let
$\AZ$ be a surjective CA and let $\BZ$ be such that $B=A \cup \{b\}$ (with $b\notin A$) and $\forall x\in \bz, G(x)=F(x')$ where $x'$ is obtained by substituting every occurrence of $b$ in $x$ with the symbol $a\in A$. Then $G$ is not surjective and $G(\bz)=\az$.

Decidability questions about properties of the limit sets received great attention. One of the most important undecidability results, due to Kari, is the following one.

\begin{theorem} \label{nilp}  \cite{Kari92} Nilpotency is undecidable for CA.
\end{theorem}

Nilpotency remains undecidable also under the additional condition of a spreading state.  Nilpotency is the basis to prove the undecidability of most of the undecidable properties of CA. In particular, Kari showed that (the problem to decide) nilpotency is the easiest problem among all decision problems on the limit sets. 

\begin{theorem}  \label{rice} \cite{Kari94} Every nontrivial property of CA limit sets is undecidable.
\end{theorem}

For example, by Theorem \ref{rice}, every nontrivial property which regards the language
$\mathcal{L}(\Omega_F)$ is undecidable. Kari's Theorem  does not concern properties of the dynamics of CA on the limit set. Here we investigate decidability questions about properties of  \emph{limit set dynamical systems} or properties of \emph{limit set dynamics}. 

\begin{definition} A property $\mathcal{P}$ is a \emph{property of limit set dynamics} if and only if  the following condition holds: if  $\#\AZ\in\mathcal{P}$ and $\BZ$ is a cellular automaton such that
$\Omega_F=\Omega_G$ and $F\rst{\Omega_F}=G\rst{\Omega_G}$ then $\#\BZ\in\mathcal{P}$. 
\end{definition}
Note that, by definition, properties of limit sets are properties of limit set dynamics while the converse is not true.  It is evident that we cannot have the equivalent of Theorem \ref{rice} for limit set dynamics.  In fact, it is easy to find nontrivial properties of the limit set dynamical systems which are decidable. Consider, for example, the set of decidable properties $\mathcal{P}_n=\{\#\AZ \mid \exists x\in\az, F^n(x)=x \}$. Since every $F$-periodic point is contained in the limit set, all $\mathcal{P}_n$ are properties of the limit set dynamics. 

In the following section we will show that there is a large class of undecidable properties of the limit set dynamics. In particular, our main result concerns properties of stable CA and properties of CA which have a unique subshift attractor. To conclude this section we show some results related to
these properties.

\begin{proposition} \label{undec1} It is undecidable whether a cellular automaton has a unique subshift attractor.
\end{proposition}
\proof{Let $\AZ$ be a cellular automaton with a spreading state $s\in A$. The clopen set
$[s]$ is $F$-invariant, spreading and $\omega([s])=\{...ssss...\}$. If we could decide 
whether a cellular automaton has a unique subshift attractor the we could decide if $\omega([s])$
is the unique attractor of $\AZ$ and then we could decide if $\AZ$ is nilpotent or not.\qed
}

\begin{proposition} \label{undec2} \cite{Kari94} It is undecidable whether a cellular automaton is stable.
\end{proposition}
\proof{Assume that we can decide if some cellular automaton is stable or not. We show that it is possible to decide nilpotency. Let $\AZ$ be cellular automaton. If $\AZ$ is not stable then it is not nilpotent. If it is stable then there exist some $n\in\n$ such that $F^n(\az)=\Omega_F$. Then is sufficient to compute
all forward images of $\az$ until we reach the limit set $\Omega_F$ and then check if it is a singleton.\qed
}
 
It is open the question whether stability is a property of limit sets (i.e. it is unknown whether there is a subshift which can be limit set both of  a stable and of an unstable cellular automaton). We can show that stability is a property of limit set dynamics. This implies that, even if there exists a subshift which is both the limit set of some stable and of some unstable CA, then the dynamics of such automata on their limit sets must be distinct. 

\begin{proposition} \label{stabprop} Let $\AZ$ be a cellular automaton. Assume that there is a cellular automaton $\BZ$ such that $\Omega_F=\Omega_G$ and $F\rst{\Omega_F}=G\rst{\Omega_F}$. Then $\BZ$ is stable if and only if $\AZ$ is stable.
\end{proposition}
\proof{Let $r$ be the maximum between the radius of $\AZ$ and the radius of $\BZ$. By Proposition \ref{attr} and Proposition \ref{spreading}, there is  a clopen, $F$-invariant spreading set $U \subseteq\mathcal{B}_{2^{-r}}(\Omega_F)$ such that $\omega_F(U)=\Omega_F$. By Proposition \ref{subshift}, there is a mixing SFT
$\Sigma \subset U\subseteq \mathcal{B}_{2^{-r}}(\Omega_F)$ such that $F(\Sigma)\subseteq \Sigma$ and $\omega_F(\Sigma)=\Omega_F$. We show that $\mathcal{L}_{2r+1}(\Sigma)=\mathcal{L}_{2r+1}(\Omega_F)$.
Since $\Omega_F\subseteq \Sigma$ it is clear that $\mathcal{L}_{2r+1}(\Sigma)\supseteq\mathcal{L}_{2r+1}(\Omega_F)$. Assume that $\mathcal{L}_{2r+1}(\Sigma)\not\subseteq\mathcal{L}_{2r+1}(\Omega_F)$. Then there must be a configuration $x \in \Sigma$ such that
$x_{[-r,r]}\notin \mathcal{L}_{2r+1}(\Omega_F)$ but this would imply that $x \notin \mathcal{B}_{2^{-r}}(\Omega_F)$ which is a contradiction. Then we have $\mathcal{L}_{2r+1}(\Sigma)=\mathcal{L}_{2r+1}(\Omega_F)=\mathcal{L}_{2r+1}(\Omega_G)$ which implies that $F\rst{\Sigma}=G\rst{\Sigma}$ and $\Sigma \subseteq \bz$.
To conclude the proof it is sufficient to show that there exist $n,m\in\n$ such that $F^n(\az)\subseteq\Sigma$ and $G^m(\az) \subseteq \Sigma$ which would imply that $\AZ$ and $\BZ$ are both stable if and only if $(\Sigma,F)$ is stable and are both unstable otherwise. 

Let $W$ be the clopen set as defined in Proposition \ref{subshift}. Since $\omega_F(\az)=\Omega_F$, by compactness, we have that for every $x\in\az$ there exists $n_x\in\n$ such that $F^{n_x}(x)\in W$. Then, since $W$ is $F$-invariant, by Lemma \ref{attr}, there exists $n\in\n$ such that $F^n(\az)\subseteq W$. Then, for every $i\in\z$, it must be $F^n(\sigma^i(x))\in W$ which implies that $F^n(x)\in \Sigma$. We obtained that there exists $n\in\n$ such that $F^n(\az)\subseteq \Sigma$. 
 
Let $t$ be the order of $\Sigma$ and let $k\in\n$ be such that $2k+1\geq t$. By using the same argument above, we can show that there exists a mixing SFT $\Sigma' \subset \mathcal{B}_{2^{-k}}(\Omega_G)$
such that $\omega_G(\Sigma')=\Omega_G$ and such that $F^m(\bz)\subseteq \Sigma'$ for some $m\in\n$. We just need to show that $\Sigma' \subseteq \Sigma$ to obtain that there exists $m\in\n$ such that $F^m(\bz) \subseteq \Sigma$. Denote with $\Omega_{2k+1}$ the 
SFT approximation of order $2k+1$ of $\Omega_G$. By using the same argument above, we have
$\mathcal{L}_{2k+1}(\Sigma')=\mathcal{L}_{2k+1}(\Omega_G)=\mathcal{L}_{2k+1}(\Omega_{2k+1})$ so $\Sigma'$ is of order $t'\geq 2k+1$.
Since $\Sigma$ is of order $t\leq 2k+1$ and $\mathcal{L}_{t}(\Sigma)\supseteq\mathcal{L}_t(\Omega_G)$ it follows that $\Sigma' \subseteq \Omega_{2k+1} \subseteq \Sigma$. 
\qed
}

\section{Undecidable properties of limit set dynamics}\label{results}

In this section we show that there is a large class of properties of the dynamics on the limit set which are not decidable. In particular we show that are undecidable all properties of the limit set dynamics which are properties of stable CA or properties of CA with a unique subshift attractor. 

\begin{definition} Let $\mathcal{P}$ be a property of CA. We say that $\mathcal{P}$ is a \emph{stable property}
if $\forall \#\AZ \in \mathcal{P}$, $\AZ$ is stable.
\end{definition}

\begin{definition} Let $\mathcal{P}$ be a property of CA. We say that $\mathcal{P}$ is a \emph{unique subshift attractor property} if $\forall \#\AZ \in \mathcal{P}$, $\Omega_F$ is the unique subshift attractor of $\AZ$.
\end{definition}

To prove the undecidability of stable properties of limit set dynamics we need a preliminary result. We don't know if a not nilpotent cellular automaton with a spreading state must be unstable. Anyway, by a simple construction, given a cellular automaton
with a spreading state, we can build a new cellular automaton with a spreading state which is nilpotent (then stable) if and only if the old one is nilpotent and it is unstable otherwise.

\begin{lemma} \label{unstable} Let $\AZ$ be a CA with a spreading state. Then it is possible to construct a CA $\BZ$ with a spreading state such that $\BZ$ is nilpotent if and only if $\AZ$ is nilpotent and $\BZ$ is unstable otherwise.
\end{lemma}
\proof{Let $s\in A$ and $r\in \n$ be the spreading state and the radius of $\AZ$, respectively. Define $B=A\cup \{s'\}$ where $s'\notin A$. We define the local rule of $\BZ$ in the following way
\begin{displaymath}
g(x_1, ..., x_{2r+1}) = \left\{ \begin{array}{ll} f(x_1, ..., x_{2r+1}) & if\;\; \forall i, x_i\in A \;\;and\;\; \exists x_i \neq s  \\ 
                                                                             s' & otherwise\\
                                      \end{array} \right.
\end{displaymath}
Note that the new state $s'$ is spreading for $\BZ$ and that the only block in $A^{2r+1}$ which is mapped to $s'$ is $s^{2r+1}$. Now, it is clear that $\AZ$ is nilpotent if and only if $\BZ$ is nilpotent. Assume that $\AZ$ is not nilpotent. By compactness, it is possible to prove that there exists a configuration $x\in\az$ such that $\forall i \in \n, \forall j \in \z, F^i(x)_j \neq s$. Define the configuration $y\in \bz$ in the following way: 
 $y_{(-\infty,-1]}=x_{(-\infty,-1]}, y_{[1,\infty)}=x_{[1,\infty)}$ and $y_0=s'$. We have that $F^{-1}(y)=\emptyset, \omega(y)=\{... s's's' ...\}$ and $\forall i \in \n, \forall j \in \z, F^i(y)_j \neq s$.
For $n\in \n$ consider $z \in F^{-n}(F^n(y))$. Since $s$ is spreading in $\az$ and since $s^{2r+1}$ is the unique block in $A^{2r+1}$ which is mapped to $s'$, the only possibility is that $z_0=s'$. Moreover it is easy to check that $\forall j\in \z\setminus \{0\}, s'\neq z_j \neq s$ and $F^{-1}(z)=\emptyset$. Then $\forall n\in \n, F^n(y)\notin \Omega_G$ which implies that $\BZ$ is unstable.\qed 
}

Note that, by Proposition \ref{undec2} and Proposition \ref{stabprop}, if a property $\mathcal{P}$ is a property of all stable CA then $\mathcal{P}$ is undecidable.

\begin{theorem} \label{stable} Every nonempty stable property of limit set dynamics is undecidable.
\end{theorem}
\proof{The proof is by reduction from nilpotency. Assume that $\mathcal{P}$ is some nonempty stable property of limit set dynamics. Let
$\#\AZ\in\mathcal{P}$ and let $\BZ$ be a cellular automaton with a spreading state $s\in B$. By Lemma \ref{unstable}, we can assume that 
 $\BZ$ is stable if and only if it is nilpotent.
 
We show how to build a new cellular automaton $\CZ$ such that $\#\CZ\in \mathcal{P}$ if and only if $\BZ$ is nilpotent. 
We can build $\CZ$ by simply taking the product of $\AZ$ with $\BZ$. In detail, consider the product cellular automaton $(\az\times\bz,F\times G)$. To obtain $\CZ$ it is sufficient to recode the alphabet of $\az\times\bz$ in the following way
\begin{displaymath}
\forall a\in A, \forall b\in B, (a,b) = \left\{ \begin{array}{ll}a & \;\;if\; b=s \\ a_b &\;\;otherwise\end{array}\right.
\end{displaymath}
Since there is a $1$-to-$1$ mapping between $\cz$ and $\az\times\bz$, the local rule of $H$ on $\cz$ is naturally induced by the local rule of $F\times G$ on $\az\times\bz$. 

Now, it is not difficult to see that $\Omega_F=\Omega_H, F\rst{\Omega_F}=H\rst{\Omega_H}$ if and only if $\BZ$ is nilpotent and in this case $\#\CZ\in \mathcal{P}$.
On the contrary $\BZ$ is unstable then $\CZ$ is also unstable and $\#\CZ\notin \mathcal{P}$.\qed
}

The construction of Theorem  \ref{stable} can be used also for the unique subshift attractor case.
Also in this case note that, by Proposition \ref{undec1}, if a property $\mathcal{P}$ is a property of all CA with a unique subshift attractor then $\mathcal{P}$ is undecidable.

\begin{theorem} \label{unique} Every nonempty unique subshift attractor property of limit set dynamics is undecidable.
\end{theorem}
\proof{The proof is by reduction from nilpotency. Let $\mathcal{P}$ be some nonempty unique subshift attractor property of limit set dynamics. Let
$\#\AZ\in\mathcal{P}$ and let $\BZ$ be a cellular automaton with a spreading state $s\in B$. By using the construction of Theorem \ref{stable}, we can build a cellular automaton $\CZ$ by taking the product of $\AZ$ with $\BZ$. We show that $\#\CZ\in \mathcal{P}$ if and only if $\BZ$ is nilpotent. As shown in Theorem \ref{stable}, we have that $\Omega_F=\Omega_H$ and $F\rst{\Omega_F}=H\rst{\Omega_H}$ if and only if $\BZ$ is nilpotent. Moreover, by construction, $A \subset C$ and the clopen set $U=\{[a] \mid a\in A\}$ is $H$-invariant, spreading and $\omega_H(U)=\Omega_F$. Then if $\BZ$ is nilpotent we have that
$\omega_H(U)=\Omega_F=\Omega_H$ and $\#\CZ\in \mathcal{P}$. Otherwise $\omega_H(U)=\Omega_F\neq\Omega_H$ and $\CZ$ has two distinct subshift attractors, $\Omega_F$ and
$\Omega_H$, then $\#\CZ\notin \mathcal{P}$.\qed
}

To conclude we show some properties of limit set dynamics which are undecidable. We need the following proposition.

\begin{theorem}\label{closing} Let $\AZ$ be a cellular automaton. If $F:\Omega_F\rightarrow \Omega_F$ is closing then $\Omega_F$ is a mixing SFT.
\end{theorem}
\proof{Since $F$ is closing on $\Omega_F$, $\exists n> 0$ such that $\forall x,y \in \Omega_F, \forall i \in \z$ if $x_{[i,i+n)}=y_{[i,i+n)}$ and  $F(x)_{[i,i+2n]}=F(y)_{[i,i+2n]}$ then $x_{i+n}=y_{i+n}$.
Consider the subshift $S = \{(x,y)\mid F(x)=y\} \subseteq \Omega_F\times\Omega_F$. Let $m=\max\{n,r\}$ where $r$ is the radius of $\AZ$. Let $S_{2m+1}$ be the SFT approximation of order $2m+1$ of $S$. Consider the two projections of $S_{2m+1}$:
\begin{itemize}
\item $S'_{2m+1}=\{x \mid \exists (x,y)\in S_{2m+1}\}$
\item $S''_{2m+1}=\{y \mid \exists (x,y)\in S_{2m+1}\}$
\end{itemize}
Since $m\geq r$, we have $F(S'_{2m+1})=S''_{2m+1}$ and $\Omega_F \subseteq S'_{2m+1}$. We show that $S'_{2m+1}$ is a SFT and that $F$ restricted to 
$S'_{2m+1}$ is closing. Since $F(S'_{2m+1})=S''_{2m+1}$, it follows that $(S'_{2m+1},\sigma)$ is conjugated to $(S_{2m+1},\sigma)$ then $S'_{2m+1}$ is a SFT. Assume for absurd that there are two sequences $x,y\in S'_{2m+1}$ such that $x_n \neq y_n$, $x_{(-\infty,n)}=y_{(-\infty,n)}$ and $F(x)=F(y)$.  Then, since $m\geq n$ and $F$ is closing on $\Omega_F$ it follows that must be $x_n=y_n$ contradicting the assumption. 

Let $k$ be the order of $S'_{2m+1}$ and let $t\in\n$ such that $2t+1\geq k$. By Proposition \ref{attr}, there exists an $F$-invariant clopen set $U\subseteq\mathcal{B}_{2^{-t}}(\Omega_F)$ such that $\omega(U)=\Omega_F$. Moreover, by Proposition \ref{subshift}, $U$ contains a mixing SFT $\Sigma$ such that $F(\Sigma)\subseteq \Sigma$ and $\omega(\Sigma)=\Omega_F$. Moreover, since $2t+1$ is larger than the order of $S'_{2m+1}$, we have also $\Omega_F \subseteq \Sigma \subseteq S'_{2m+1}$. 
Now, since $F$ is closing on $S'_{2n+1}$, it follows that $F$ must be closing on $\Sigma$. Then
since $\Sigma$ is mixing, $F(\Sigma)\subseteq \Sigma$ and $F$ is closing on $\Sigma$ it follows that $F(\Sigma)=\Sigma$ which implies that $\Sigma\subseteq \Omega_F$ and then $\Omega_F=\Sigma$.\qed
}

From Theorem \ref{stable} and Theorem \ref{unique} we can easily derive the following corollary. 
\newpage
\begin{corollary} There is no algorithm that, given $\AZ$, can decide if 
\begin{itemize}
\item[$1.$] $F:\Omega_F \rightarrow \Omega_F$ is transitive.
\item[$2.$] $F:\Omega_F \rightarrow \Omega_F$ is closing,
\item[$3.$] $F:\Omega_F \rightarrow \Omega_F$ is injective,
\item[$4.$] $F:\Omega_F \rightarrow \Omega_F$ is the identity map,
\item[$5.$] $F:\Omega_F \rightarrow \Omega_F$ is expansive,
\item[$6.$] $F:\Omega_F \rightarrow \Omega_F$ is positively expansive,
\end{itemize}
\end{corollary}
\proof{By Theorem \ref{stable} and Theorem \ref{unique} it is sufficient to show that properties $1,..,6$ imply that $\AZ$ is stable or that it has a unique subshift attractor.
\begin{itemize}
\item[$1.$] If $F:\Omega_F \rightarrow \Omega_F$ is transitive then $\Omega_F$ is the unique attractor of $\AZ$ and, in particular, it is the unique subshift attractor.
\item[$2.$] If $F:\Omega_F \rightarrow \Omega_F$ is closing then, by Theorem \ref{closing}, $\Omega_F$ is a mixing SFT then  $\AZ$ must be stable.
\item[$3.$] If $F:\Omega_F \rightarrow \Omega_F$ is injective then, since $F$ is surjective on $\Omega_F$, it must be invertible and then closing. 
\item[$4.$] If $F:\Omega_F \rightarrow \Omega_F$ is the identity map then $F$ must be injective on $\Omega_F$.
\item[$5.$] If $F:\Omega_F \rightarrow \Omega_F$ is expansive then $F$ must be injective on
$\Omega_F$ and then closing and transitive.
\item[$6.$] If $F:\Omega_F \rightarrow \Omega_F$ is positively expansive then $F$ must be closing on
$\Omega_F$ and transitive.
\end{itemize}
}

\section{Concluding remarks}\label{conclusions}
In this paper we proved that any property of limit set dynamics is undecidable, if it implies stability or the existence of a unique subshift attractor. As examples of properties which imply stability we have closing (which implies that the limit set is a mixing SFT), injectivity, expansivity, positively expansiveness and identity (all of which imply closing). As examples of properties which imply the existence of a unique subshift attractor we have transitivity, expansivity and positively expansiveness (expansive and positively expansive endomorphisms of mixing SFTs are transitive).
From Theorem \ref{stable} and \ref{unique} we can conclude that all such properties are undecidable. We remark that, since surjectivity is not a property of limit set dynamics (and it is decidable), if we restrict to only surjective CA then we cannot derive any conclusion from our theorems. In particular we cannot conclude anything about the decidability of transitivity, expansivity and positively expansiveness (it is already known that closing, injectivity and identity are decidable for surjective CA).

Our main undecidability proofs are by reduction from nilpotency. Note that a nilpotent CA is stable and it has a unique subshift attractor. Then (the problem to decide) nilpotency is the easiest problem among all decision problems on the limit set dynamics of stable CA and of CA with a unique subshift attractor.

We conclude the paper by raising a question. It is not clear how stability is related to the existence of a unique subshift attractor. To our knowledge there are no examples of stable CA with two distinct subshift attractors. For a wide class of stable CA it is possible to prove that they have a unique subshift attractor (in particular surjective CA, see \cite{FK07}) but the general question is open. If stable CA have a unique subshift attractor then Lemma \ref{unstable} would be useless and we could derive Theorem \ref{stable} as a corollary of Theorem \ref{unique}.

\begin{question} Is there any stable CA with two distinct subshift attractors?
\end{question}


\newpage\null
\end{document}